\input amstex
\documentstyle{amsppt}
%
\catcode`@=11
\redefine\output@{%
  \def\break{\penalty-\@M}\let\par\endgraf
  \ifodd\pageno\global\hoffset=105pt\else\global\hoffset=8pt\fi  
  \shipout\vbox{%
    \ifplain@
      \let\makeheadline\relax \let\makefootline\relax
    \else
      \iffirstpage@ \global\firstpage@false
        \let\rightheadline\frheadline
        \let\leftheadline\flheadline
      \else
        \ifrunheads@ 
        \else \let\makeheadline\relax
        \fi
      \fi
    \fi
    \makeheadline \pagebody \makefootline}%
  \advancepageno \ifnum\outputpenalty>-\@MM\else\dosupereject\fi
}
\def\Beta{\mathchar"0\hexnumber@\rmfam 42}
\redefine\mm@{2010} 
\catcode`\@=\active
\nopagenumbers
\chardef\textvolna='176

\chardef\bigalpha='013
\def\negskp{\hskip -2pt}

\chardef\degree="5E
\def\compos{\,\raise 1pt\hbox{$\sssize\circ$} \,}

\def\vtrule{\vrule height 12pt depth 6pt}
\def\blue#1{#1}

\gdef\darkred#1{#1}
\catcode`#=11\def\diez{#}\catcode`#=6
\catcode`&=11\catcode`&=4
\catcode`_=11\catcode`_=8
\catcode`\^=11\catcode`\^=7
\catcode`~=11\catcode`~=\active
\catcode`\%=11\def\procent{
\def\mycite#1{\cite{\blue{#1}}\immediate\special{ps:
     ShrHPSdict begin /ShrBORDERthickness 0 def}}
\def\myciterange#1#2#3#4{\cite{\blue{#2#3#4}}\immediate\special{ps:
     ShrHPSdict begin /ShrBORDERthickness 0 def}}
\def\mytag#1{%
    \tag#1}
\def\mythetag#1{\thetag{\blue{#1}}\immediate\special{ps:
     ShrHPSdict begin /ShrBORDERthickness 0 def}}
\def\myrefno#1{\no#1}
\def\myhref#1#2{\blue{#2}\immediate\special{ps:
     ShrHPSdict begin /ShrBORDERthickness 0 def}}
\def\myEarXivlink{\myhref{http://arXiv.org}{http:/\negskp/arXiv.org}}
\def\myGeoCities{\myhref{http://www.geocities.com}{GeoCities}}
\def\mytheorem#1{\csname proclaim\endcsname{Theorem #1}}
\def\mytheoremwithtitle#1#2{\csname proclaim\endcsname{Theorem #1#2}}
\def\mythetheorem#1{\blue{#1}\immediate\special{ps:
     ShrHPSdict begin /ShrBORDERthickness 0 def}}
\def\mylemma#1{\csname proclaim\endcsname{Lemma #1}}
\def\mylemmawithtitle#1#2{\csname proclaim\endcsname{Lemma #1#2}}
\def\mythelemma#1{\blue{#1}\immediate\special{ps:
     ShrHPSdict begin /ShrBORDERthickness 0 def}}
\def\mycorollary#1{\csname proclaim\endcsname{Corollary #1}}
\def\mythecorollary#1{\blue{#1}\immediate\special{ps:
     ShrHPSdict begin /ShrBORDERthickness 0 def}}
\def\mydefinition#1{\definition{Definition #1}}
\def\mythedefinition#1{\blue{#1}\immediate\special{ps:
     ShrHPSdict begin /ShrBORDERthickness 0 def}}
\def\myconjecture#1{\csname proclaim\endcsname{Conjecture #1}}
\def\myconjecturewithtitle#1#2{\csname proclaim\endcsname{Conjecture #1#2}}
\def\mytheconjecture#1{\blue{#1}\immediate\special{ps:
     ShrHPSdict begin /ShrBORDERthickness 0 def}}
\def\myproblem#1{\csname proclaim\endcsname{Problem #1}}
\def\myproblemwithtitle#1#2{\csname proclaim\endcsname{Problem #1#2}}
\def\mytheproblem#1{\blue{#1}\immediate\special{ps:
     ShrHPSdict begin /ShrBORDERthickness 0 def}}
\def\mytable#1{Table #1}
\def\mythetable#1{\blue{#1}\immediate\special{ps:
     ShrHPSdict begin /ShrBORDERthickness 0 def}}
\def\myanchortext#1#2{#2}
\def\mytheanchortext#1#2{\blue{#2}\immediate\special{ps:
     ShrHPSdict begin /ShrBORDERthickness 0 def}}
\font\eightcyr=wncyr8
\pagewidth{360pt}
\pageheight{606pt}
\topmatter
\title
On a simplified version of Hadamard's maximal 
determinant problem\endtitle
\rightheadtext{On a simplified version of Hadamard's problem}
\author
Ruslan Sharipov
\endauthor
\address Bashkir State University, 32 Zaki Validi street, 450074 Ufa, Russia
\endaddress
\email
\myhref{mailto:r-sharipov\@mail.ru}{r-sharipov\@mail.ru}
\endemail
\abstract
Hadamard's maximal determinant problem consists in finding the maximal value
of the determinant of a square $n\times n$ matrix whose entries are plus or 
minus ones. This is a difficult mathematical problem which is not yet solved.
In the present paper a simplified version of the problem is considered and
studied numerically.
\endabstract
\endtopmatter
\loadbold
\TagsOnRight
\document

\head
1. Introduction.
\endhead
     Hadamard's maximal determinant problem was first published in \mycite{1}
by Jacques Salomon Hadamard\footnotemark\ in 1893. Let's denote
$$
\hskip -2em
M(n)=\max\limits_{a_{ij}\,=\,\pm 1}\vmatrix a_{11} & \hdots & a_{1n}\\
\vdots & \ddots & \vdots\\
a_{n1} & \hdots & a_{nn}\endvmatrix.
\mytag{1.1}
$$
\footnotetext{\,Mark Kac in \mycite{2} supposed that the problem belongs 
to Maurice Ren\'e Fr\'echet. Fr\'echet was a student of Hadamard at secondary 
school Lyc\'ee Buffon in Paris, so the version has some ground. This version
is supported in \mycite{3}. With the reference to \mycite{3} in the Russian
segment of Wikipedia the problem is called Fr\'echet's problem of maximal 
determinant (see \mycite{4}).}\adjustfootnotemark{-1}\noindent
The problem is to find $M(n)$ in \mythetag{1.1} for each $n\in\Bbb N$. There is 
a closely related problem of finding $d_n$, where
$$
\hskip -2em
d_n=\max\limits_{a_{ij}\in\,\{0,1\}}\vmatrix a_{11} & \hdots & a_{1n}\\
\vdots & \ddots & \vdots\\
a_{n1} & \hdots & a_{nn}\endvmatrix.
\mytag{1.2}
$$
According to \mycite{5}, the relation of two problems \mythetag{1.1} and \mythetag{1.2} 
is given by the formula
$$
\hskip -2em
M(n)=2^{n-1}\,D(n)\text{, \ where \ }D(n)=d_{n-1}.
\mytag{1.3}
$$
There are various estimates for $M(n)$, $D(n)$, and $d_n$ (see \mycite{1}, \mycite{5},
and \myciterange{5}{6}{--}{12}). However the exact values of $M(n)$ and $D(n)$ in
\mythetag{1.3} are known only for $n\leqslant 21$. The exact values of $d_n$ are known 
for $n\leqslant 20$ (see \mycite{5}).\par
    The goal of the present paper is neither to improve existing estimates not to give 
new ones. Here we consider a different problem which share some features of Hadamard's 
problem, but is not dependent of it. Acting like in \mythetag{1.2}, let's denote through 
$A_n$ a square $n\times n$ matrix whose elements are zeros and ones:
$$
\hskip -2em
A_n=\Vmatrix a_{11} & \hdots & a_{1n}\\
\vdots & \ddots & \vdots\\
a_{n1} & \hdots & a_{nn}\endVmatrix\text{, \ where \ }a_{ij}\in\,\{0,1\}.
\mytag{1.4}
$$
The problem then is formulated as follows.
\myproblem{1.1} Starting from the unit matrix $A_1=\Vert 1\Vert$ for $n=1$, find 
a sequence $A_n$ of square $n\times n$ matrices of the form \mythetag{1.4} each next of 
which comprises the previous one as its upper left diagonal block and is of maximal 
determinant in the so described class.
\endproclaim
We shall call this problem the simplified Hadamard maximal determinant problem until
some more convenient name will be suggested.
\head
2. Computer code for solving the problem.
\endhead
     The first three matrices $A_1$, $A_3$, $A_3$, in the series are easily written
$$
\xalignat 3
&\hskip -2em
A_1=\Vert 1\Vert,
&&A_2=\Vmatrix 1 & 0\\ 1 & 1\endVmatrix,
&&A_3=\Vmatrix 1 & 0 & 1\\ 1 & 1 & 0\\ 0 & 1 & 1\endVmatrix.
\mytag{2.1}
\endxalignat
$$
The matrices \mythetag{2.1} provide the following series of maximal determinants
$b_i=\det A_i$: 
$$
\xalignat 3
&\hskip -2em
b_1=1, &&b_2=1, &&b_3=2.
\mytag{2.2}
\endxalignat
$$
Further terms in \mythetag{2.1} and \mythetag{2.2} are produced computationally.
Since the definition of the matrices $A_n$ in Problem~\mytheproblem{1.1} is by induction,
the matrix $A_{n+1}$ looks like
$$
\hskip -2em
A_{n+1}=\Vmatrix a_{11} & \hdots & a_{1n} & y_1\\
\vdots & \ddots & \vdots & \vdots\\
a_{n1} & \hdots & a_{nn} & y_n \\
x_1    & \hdots & x_n & z\endVmatrix.
\mytag{2.3}
$$
When computing $A_{n+1}$ its upper left diagonal block is already known, it coincides 
with $A_n$. The determinant of the matrix \mythetag{2.3} is given by the formula
$$
\hskip -2em
\det A_{n+1}=z\,\det A_n+\sum^n_{i=1}\sum^n_{j=1} b_{ij}\,x_i\,y_j.
\mytag{2.4}
$$
Since $\det A_n>0$, the choice of $z$ in the formula \mythetag{2.4} is obvious from the 
maximal determinant condition in Problem~\mytheproblem{1.1}: 
$$
\hskip -2em
z=1.
\mytag{2.5}
$$
Substituting \mythetag{2.5} into \mythetag{2.4}, we derive 
$$
\hskip -2em
\det A_{n+1}=\det A_n+\sum^n_{i=1}\sum^n_{j=1} b_{ij}\,x_i\,y_j.
\mytag{2.6}
$$
The double sums in \mythetag{2.6} constitute a bilinear form with respect to the
variables $x_1,\,\ldots,\,x_n$ and $y_1,\,\ldots,\,y_n$. The coefficients $b_{ij}$ 
of this bilinear form are expressed through the entries of the matrix $A_n$. The
following code was used for computing the determinant $\det(A_{n+1})$ in 
\mythetag{2.6}.\par
\medskip
\parshape 1 10pt 350pt 
{\tt\noindent\darkred{function\_to\_maximize(PM):=block}\newline
\darkred{\ ([n,M,i,j,D],}\newline
\darkred{\ \ n:length(PM)+1,}\newline
\darkred{\ \ M:zeromatrix(n,n),}\newline
\darkred{\ \ for i:1 step 1 thru n-1 do}\newline
\darkred{\ \ \ for j:1 step 1 thru n-1 do}\newline
\darkred{\ \ \ \ M[i,j]:PM[i,j],}\newline
\darkred{\ \ for i:1 step 1 thru n-1 do}\newline
\darkred{\ \ \ (}\newline
\darkred{\ \ \ \ M[n,i]:x[i],}\newline
\darkred{\ \ \ \ M[i,n]:y[i]}\newline
\darkred{\ \ \ )}\newline
\darkred{\ \ M[n,n]:1,}\newline
\darkred{\ \ M\_max:copy(M),}\newline
\darkred{\ \ D:determinant(M),}\newline
\darkred{\ \ D:ratexpand(D),}\newline
\darkred{\ \ return(D)}\newline
\darkred{\ )\$}}\par
\medskip
This code defines a function with the name {\tt \darkred{function\_to\_maximize}}, 
where {\tt\darkred{[n,M,}} {\tt\darkred{i,j,D]}} is the list of its local variables. 
The function {\tt \darkred{function\_to\_maximize}} uses one global variable 
{\tt\darkred{M\_max}} in order to output the matrix \mythetag{2.3} with $z=1$. The
argument {\tt\darkred{PM}} of this function is used in order to load the input matrix
$A_n$ of the form \mythetag{1.4}. The result of this function is the expression of the 
form \mythetag{2.6} that should be maximized by choosing proper values for the
variables $x_i$ and $y_i$.\par
     The above code is written in Maxima programming language. Maxima is a Computer 
Algebra System available for Linux, Windows, and MacOS. I run Maxima version 5.42.2 
on the Linux platform Ubuntu 16.04 LTS.\par
      Another code is used in order to maximize the expression \mythetag{2.6}. It is
as follows.\par
\medskip
\parshape 1 10pt 350pt 
{\tt\noindent\darkred{compute\_max\_det(n,F):=block}\newline
\darkred{\ ([i,j,P,r,rr],}\newline
\darkred{\ \ r:0,}\newline
\darkred{\ \ if n=0}\newline
\darkred{\ \ \ then }\newline
\darkred{\ \ \ \ (}\newline
\darkred{\ \ \ \ \ if F>D\_max}\newline
\darkred{\ \ \ \ \ \ then}\newline
\darkred{\ \ \ \ \ \ \ (}\newline
\darkred{\ \ \ \ \ \ \ \ D\_max:F,}\newline
\darkred{\ \ \ \ \ \ \ \ r:1}\newline
\darkred{\ \ \ \ \ \ \ )}\newline
\darkred{\ \ \ \ )}\newline
\darkred{\ \ \ else}\newline
\darkred{\ \ \ \ (}\newline
\darkred{\ \ \ \ \ for i:0 step 1 thru 1 do}\newline
\darkred{\ \ \ \ \ \ for j:0 step 1 thru 1 do}\newline
\darkred{\ \ \ \ \ \ \ (}\newline
\darkred{\ \ \ \ \ \ \ \ P:psubst([x[n]=i,y[n]=j],F),}\newline
\darkred{\ \ \ \ \ \ \ \ rr:compute\_max\_det(n-1,P),}\newline
\darkred{\ \ \ \ \ \ \ \ if rr=1}\newline
\darkred{\ \ \ \ \ \ \ \ \ then}\newline
\darkred{\ \ \ \ \ \ \ \ \ \ (}\newline
\darkred{\ \ \ \ \ \ \ \ \ \ \ x\_max[n]:i,}\newline
\darkred{\ \ \ \ \ \ \ \ \ \ \ y\_max[n]:j,}\newline
\darkred{\ \ \ \ \ \ \ \ \ \ \ r:1}\newline
\darkred{\ \ \ \ \ \ \ \ \ \ )}\newline
\darkred{\ \ \ \ \ \ \ )}\newline
\darkred{\ \ \ \ ),}\newline
\darkred{\ \ return(r)}\newline
\darkred{\ )\$}}\par
\medskip
\noindent This code defines the function {\tt\darkred{compute\_max\_det}} with 
local variables {\tt\darkred{i,j,P,r,rr}} and with two arguments {\tt\darkred{n,F}}.
The argument {\tt\darkred{F}} is used in order to load the expression \mythetag{2.6},
the argument {\tt\darkred{n}} admits the value of the number $n$ in \mythetag{2.6}.
The function {\tt\darkred{compute\_max\_det}} uses three global variables 
{\tt\darkred{D\_max,x\_max,y\_max}}. The variable {\tt\darkred{D\_max}} is used in
order to output the maximal value of expression \mythetag{2.6}, i\.e\. the maximal 
value of the determinant. The global variables {\tt\darkred{x\_max,y\_max}} are used
as undeclared arrays. Through them we output the values of $x_1,\,\ldots,\,x_n$ and
$y_1,\,\ldots,\,y_n$ at which the maximum of the determinant \mythetag{2.6} is attained.
Substituting these values along with \mythetag{2.5} into \mythetag{2.4}, we get the matrix
$A_{n+1}$ at which the maximum of the determinant is reached.\par
     The main code for computing several matrices $A_n$ looks like a loop with respect
to the integer variable $n$. We start with the matrix $A_3$ in \mythetag{2.1}.
\medskip
\parshape 1 10pt 350pt 
{\tt\noindent\darkred{file\_output\_append:true\$}\newline
\darkred{MM:matrix([1,0,1],[1,1,0],[0,1,1])\$}\newline
\darkred{for n:3 step 1 thru 14 do}\newline
\darkred{\ (}\newline
\darkred{\ \ FF:function\_to\_maximize(MM),}\newline
\darkred{\ \ N\_max:length(listofvars(FF))/2,}\newline
\darkred{\ \ D\_max:0,}\newline
\darkred{\ \ compute\_max\_det(N\_max,FF),}\newline
\darkred{\ \ MM:copy(M\_max),}\newline
\darkred{\ \ for i:1 step 1 thru N\_max do }\newline
\darkred{\ \ \ MM:psubst([x[i]=x\_max[i],
             y[i]=y\_max[i]],MM),}\newline
\darkred{\ \ stringout("output\_file",[n+1,D\_max,MM])}\newline
\darkred{\ )\$}}\par
\medskip
\noindent The matrices $A_n$ along with their determinants are written to the 
{\tt\darkred{output\_file}}\,.\par
\head
3. The result of computations and conclusions.
\endhead
     The above code was run once in a loop for $n$ from $n=3$ through $n=14$. As 
a result the matrices $A_4$, $A_5$, $A_6$, $A_7$, $A_8$, $A_9$, $A_{10}$, $A_{11}$, 
$A_{12}$, $A_{13}$, $A_{14}$, $A_{15}$ were computed in addition to the matrices 
\mythetag{2.1}. According to the statement of Problem~\mytheproblem{1.1}, these
matrices are enclosed in each other like "matryoshkas" (nesting dolls) in the form 
of upper left diagonal blocks.  Therefore it is sufficient to typeset only the last
one of them, i\.\,e\. the matrix $A_{15}$: 
$$
\hskip -2em
A_{15}=\Vmatrix 
1&0&1&1&0&0&0&0&0&0&1&0&0&1&1\\
1&1&0&0&0&1&0&1&1&1&1&0&0&0&0\\
0&1&1&0&1&0&0&1&0&0&1&1&0&1&0\\
0&1&0&1&1&1&0&0&0&0&0&0&1&0&1\\
1&0&0&0&1&1&1&1&0&0&0&0&0&1&0\\
0&0&1&0&0&1&1&0&1&0&0&0&1&1&0\\
0&1&0&1&0&0&1&1&0&0&1&0&1&1&0\\
0&0&1&1&0&1&0&1&0&1&0&1&1&0&0\\
0&0&0&1&1&0&0&1&1&1&0&0&0&1&1\\
1&1&1&0&1&0&1&0&0&1&0&0&1&0&0\\
0&0&0&0&1&1&1&0&0&1&1&1&0&0&1\\
1&1&0&1&0&0&1&0&1&0&0&1&0&0&0\\
1&0&0&0&1&0&0&1&1&0&1&1&1&0&1\\
1&1&0&0&0&1&0&0&0&1&0&1&1&1&1\\
0&1&1&0&0&0&1&1&0&0&0&0&0&0&1
\endVmatrix.
$$\par
     Along with the matrices we obtain their determinants which are 
maximal in their classes. Indeed, according to the statement of 
Problem~\mytheproblem{1.1}, we have 
$$
\hskip -2em
b_{n+1}=\det(A_{n+1})=\max\Sb x_1,\,\ldots,\,x_n\in\{0,1\}\\
y_1,\,\ldots,\,y_n\in\{0,1\}\endSb\ 
\vmatrix a_{11} & \hdots & a_{1n} & y_1\\
\vdots & \ddots & \vdots & \vdots\\
a_{n1} & \hdots & a_{nn} & y_n \\
x_1    & \hdots & x_n & 1\endvmatrix.
\mytag{3.1}
$$
For the purposes of comparison with the original Hadamard's problem we provide 
the quantities $b_n$ along with the quantities $d_n$ from \mythetag{1.3}:
$$
\vcenter{\vbox{\hsize 15cm
\offinterlineskip\settabs\+\indent
\vtrule
\hss\ $d_n$\hss&\vtrule
\hss\ \ 1\ \ \ \hss&\vtrule
\hss\ \ 1\ \ \ \hss&\vtrule
\hss\ \ 2\ \ \ \hss&\vtrule
\hss\ \ 3\ \ \ \hss&\vtrule
\hss\ \ 5\ \ \ \hss&\vtrule
\hss\ \ 9\ \ \ \hss&\vtrule
\hss\ \ 32\ \ \ \hss&\vtrule
\hss\ \ 56\ \ \ \hss&\vtrule
\hss\ \ 144\ \ \ \hss&\vtrule
\cr\hrule 
\+\vtrule
\hss\,\ $n$\hss&\vtrule
\hss\ \ 1\hss&\vtrule
\hss\ \ 2\hss&\vtrule
\hss\ \ 3\hss&\vtrule
\hss\ \ 4\hss&\vtrule
\hss\ \ 5\hss&\vtrule
\hss\ \ 6\hss&\vtrule
\hss\ \ 7\hss&\vtrule
\hss\ \ 8\hss&\vtrule
\hss\ \ 9\hss&\vtrule
\cr\hrule
\+\vtrule
\hss\ $b_n$\hss&\vtrule
\hss\ \ 1\hss&\vtrule
\hss\ \ 1\hss&\vtrule
\hss\ \ 2\hss&\vtrule
\hss\ \ 3\hss&\vtrule
\hss\ \ 5\hss&\vtrule
\hss\ \ 9\hss&\vtrule
\hss\ \ 18\hss&\vtrule
\hss\ \ 40\hss&\vtrule
\hss\ \ 96\hss&\vtrule
\cr\hrule
\+\vtrule
\hss\ $d_n$\hss&\vtrule
\hss\ \ 1\hss&\vtrule
\hss\ \ 1\hss&\vtrule
\hss\ \ 2\hss&\vtrule
\hss\ \ 3\hss&\vtrule
\hss\ \ 5\hss&\vtrule
\hss\ \ 9\hss&\vtrule
\hss\ \ 32\hss&\vtrule
\hss\ \ 56\hss&\vtrule
\hss\ \ 144\hss&\vtrule
\cr\hrule
}}
$$
\ \ 
$$
\vcenter{\vbox{\hsize 12cm
\offinterlineskip\settabs\+\indent
\vtrule
\hss\ $d_n$\hss&\vtrule
\hss\ \ 320\ \ \hss&\vtrule
\hss\ \ 1458\ \ \hss&\vtrule
\hss\ \ 3645\ \ \hss&\vtrule
\hss\ \ 9477\ \ \hss&\vtrule
\hss\ \ 25515\ \ \hss&\vtrule
\hss\ \ 131073\ \ \ \hss&\vtrule
\cr\hrule 
\+\vtrule
\hss\,\ $n$\hss&\vtrule
\hss\ \ 10\hss&\vtrule
\hss\ \ 11\hss&\vtrule
\hss\ \ 12\hss&\vtrule
\hss\ \ 13\hss&\vtrule
\hss\ \ 14\hss&\vtrule
\hss\ \ 15\hss&\vtrule
\cr\hrule
\+\vtrule
\hss\ $b_n$\hss&\vtrule
\hss\ \ 220\hss&\vtrule
\hss\ \ 604\hss&\vtrule
\hss\ \ 1608\hss&\vtrule
\hss\ \ 4734\hss&\vtrule
\hss\ \ 14898\hss&\vtrule
\hss\ \ 45034\hss&\vtrule
\cr\hrule
\+\vtrule
\hss\ $d_n$\hss&\vtrule
\hss\ \ 320\hss&\vtrule
\hss\ \ 1458\hss&\vtrule
\hss\ \ 3645\hss&\vtrule
\hss\ \ 9477\hss&\vtrule
\hss\ \ 25515\hss&\vtrule
\hss\ \ 131073\hss&\vtrule
\cr\hrule
}}
$$
\smallskip
The inequality $b_n\leqslant d_n$ observed in the above tables is not surprising.
The quantities $d_n$ in \mythetag{1.2} are defined as total maxima with respect to 
all entries of the matrices, while $b_n$ in \mythetag{3.1} are partial maxima 
calculated for the case where entries of the upper left blocks of the matrices 
are fixed.\par
     Note that the matrix $A_2$ in the sequence \mythetag{2.1} is not unique. There 
are two other equivalent options for choosing this matrix:
$$
\xalignat 2
&\hskip -2em
A_2=\Vmatrix 1 & 1\\ 0 & 1\endVmatrix,
&&A_2=\Vmatrix 1 & 0\\ 0 & 1\endVmatrix.
\mytag{3.2}
\endxalignat
$$
These two options \mythetag{3.2} are associated with the same value of $b_2$.
Such a non-uniqueness can happen for $n>2$ as well. Therefore the solution
of the simplified Hadamard maximal determinant problem is a collection of 
sequences rather than a single sequence of matrices. Each two sequences of
this collection share some initial part thus producing a tree structure 
in the collection as a whole.\par
     The above code produces only one sequence of the collection. This means
that the research of the simplified Hadamard maximal determinant problem should 
be continued. Probably this would contribute to the solution of the original
Hadamard's maximal determinant problem as well.
\head
4. Dedicatory.
\endhead
     This paper is dedicated to my sister Svetlana Abdulovna Sharipova. 
\head
5. Acknowledgments.
\endhead
     I am grateful to Nazim Mehdiyev from Baku, the capital of Azerbaijan, who
drew my attention to Fr\'echet's problem of maximal determinant which coincides
with Hadamard's maximal determinant problem as it was explained above in the 
introductory section of this paper.\par 

\enddocument
\end